\documentclass{amsart}
\usepackage{amsmath,amssymb,color}
\usepackage{amscd,latexsym}
\usepackage[all]{xy}

\newcommand{\be}{\mathbf{1}}

\newcommand{\id}{\mathrm{id}}

\newcounter{number}[section]

\newenvironment{nummer}{\refstepcounter{number}
{\noindent\arabic{section}.\arabic{number}}}{}

\newcommand{\bn}{\begin{nummer} \rm}
\newcommand{\en}{\end{nummer}}

\newenvironment{thms}{\noindent {\sc Theorem:} \it}{}
\newenvironment{lms}{\noindent {\sc Lemma:} \it}{}
\newenvironment{props}{\noindent {\sc Proposition:} \it}{}

\newenvironment{cors}{\noindent {\sc Corollary:} \it}{}

\newenvironment{rems}{\noindent {\sc Remark:}}{}

\newenvironment{nproof}{\noindent {\sc Proof:}}{\mbox{}\hfill
\rule[-.2ex]{.25em}{1.8ex}}

\subjclass[2000]{46L35, 46L85}

\begin{document}
\author{Wilhelm Winter}

\address{School of Mathematical Sciences,  University of Nottingham, Nottingham, NG7 2RD}
\email{wilhelm.winter@nottingham.ac.uk}

\date{\today}

\thanks{{\it Supported by:} EPSRC First Grant EP/G014019/1}

\keywords{strongly self-absorbing $C^*$-algebra, Jiang--Su algebra}

\title[Strongly self-absorbing $C^{*}$-algebras are $\mathcal{Z}$-stable]{Strongly self-absorbing $C^{*}$-algebras are $\mathcal{Z}$-stable}

\begin{abstract}
We prove the title. This characterizes the Jiang--Su algebra $\mathcal{Z}$ as the uniquely determined initial object in the category of strongly self-absorbing $C^{*}$-algebras.
\end{abstract}

\maketitle

\setcounter{section}{-1}

\section{Introduction}   
\label{introduction}

\noindent
A separable unital $C^{*}$-algebra $\mathcal{D} \neq \mathbb{C}$ is called strongly self-absorbing, if there is an isomorphism $\mathcal{D} \to \mathcal{D} \otimes \mathcal{D}$ which is approximately unitarily equivalent to the first factor embedding, cf.\ \cite{TomsWin:ssa}. The interest in such algebras largely arises from Elliott's program to classify nuclear $C^{*}$-algebras by $K$-theoretic invariants. In fact, examples suggest that classification will only be possible up to $\mathcal{D}$-stability (i.e., up to tensoring with $\mathcal{D}$) for a strongly self-absorbing $\mathcal{D}$, cf.\ \cite{Toms:classproblem}, \cite{EllToms:regularity}, \cite{Win:localizingEC}. While the known strongly self-absorbing examples are quite well understood, and are entirely classified, it remains an open problem whether these are the only ones. From a more general perspective, the question is in how far abstract properties allow for comparison with concrete examples. For nuclear $C^{*}$-algebras, this question prominently manifests itself as the UCT problem (i.e., is every nuclear $C^{*}$-algebra  $KK$-equivalent to a commutative one); a positive answer even in the special setting of strongly self-absorbing $C^{*}$-algebras would be highly satisfactory, and likely shed light on the general case.      

In this note we shall be concerned with a closely related interpretation of the aforementioned question: we will show that any strongly self-absorbing $C^{*}$-algebra $\mathcal{D}$ admits a unital embedding of a specific example, the Jiang--Su algebra $\mathcal{Z}$ (we refer to \cite{JiaSu:Z} and to \cite{RorWin:Z-revisited} for an introduction and various characterizations of $\mathcal{Z}$). It then follows immediately that $\mathcal{D}$ is in fact $\mathcal{Z}$-stable. The result answers some problems left open in \cite{TomsWin:ssa} and in \cite{DadRor:ssa-projection}; in particular it implies that strongly self-absorbing $C^{*}$-algebras are always $K_{1}$-injective. Moreover, it shows that the Jiang--Su algebra is an initial object in the category of strongly self-absorbing $C^{*}$-algebras (with unital $*$-homomorphisms); there can only be one such initial object, whence $\mathcal{Z}$ is characterized this way. It is interesting to note that the Cuntz algebra $\mathcal{O}_{2}$ is the uniquely determined final object in this category, and that $\mathcal{O_{\infty}}$ can be characterized as the initial object in the category of infinite strongly self-absorbing $C^{*}$-algebras. 

The proof of our main result builds on ideas from \cite{RorWin:Z-revisited} and from \cite{DadRor:ssa-projection}, where the problem was settled in the case where $\mathcal{D}$ contains a nontrivial projection.


\section{Small elementary tensors}

\noindent
In this section, we generalize a technical result from \cite{DadRor:ssa-projection} to a setting that does not require the existence of projections, see Lemma~\ref{C} below. We refer to \cite{Ror:Z-absorbing} for a brief account of the Cuntz semigroup.

\bn
\label{A}
\begin{props}
Let $A$ be a unital $C^{*}$-algebra, $0 \le g \le \be_{A}$. 

Then, for any $0 \neq n \in \mathbb{N}$, we have
\begin{eqnarray*}
\be_{A^{\otimes n}} - g^{\otimes n} & \ge & (\be_{A} - g) \otimes g \otimes \ldots \otimes g \\
& & + g \otimes (\be_{A} - g) \otimes g \otimes \ldots \otimes g \\
&&  \vdots \\
&& + g \otimes \ldots \otimes g \otimes (\be_{A} - g).
\end{eqnarray*}
\end{props}

\begin{nproof}
The statement is trivial for $n=1$. Suppose now we have shown the assertion for some $0 \neq n \in \mathbb{N}$. We obtain
\begin{eqnarray*}
\be_{A^{\otimes (n+1)}} - g^{\otimes (n+1)} & = &  \be_{A^{\otimes n}} \otimes g - g^{\otimes n} \otimes g + \be_{A^{\otimes n}} \otimes (\be_{A} - g) \\
& = & (\be_{A^{\otimes n}}  - g^{\otimes n}) \otimes g + \be_{A^{\otimes n}} \otimes (\be_{A} - g) \\
& \ge & ((\be_{A} - g) \otimes g \otimes \ldots \otimes g) \otimes g \\
& & + (g \otimes (\be_{A} - g) \otimes g \otimes \ldots \otimes g) \otimes g \\
&&  \vdots \\
&& + (g \otimes \ldots \otimes g \otimes (\be_{A} - g)) \otimes g \\
&& + g^{\otimes n} \otimes (\be_{A} - g), 
\end{eqnarray*}
where for the inequality we have used our induction hypothesis as well as the fact that $\be_{A^{\otimes n}} \otimes (\be_{A} - g) \ge g^{\otimes n} \otimes (\be_{A} - g)$. Therefore, the statement also holds for $n+1$.
\end{nproof}
\en

\bn
\label{F}
\begin{props}
Let $\mathcal{D}$ be strongly self-absorbing,  $0 \le d \le \be_{\mathcal{D}}$.

Then, for any $0 \neq k \in \mathbb{N}$,
\[
[\be_{\mathcal{D}^{\otimes k}} - d^{\otimes k}] \le k \cdot [(\be_{\mathcal{D}}- d)  \otimes \be_{\mathcal{D}^{\otimes (k-1)}}] \mbox{ in } W(\mathcal{D}^{\otimes k}).
\]
\end{props}

\begin{nproof}
The assertion holds trivially for $k = 1$. Suppose now it has been verified for some $k \in \mathbb{N}$. Then,  
\begin{eqnarray*}
[\be_{\mathcal{D}} - d^{\otimes (k+1)}] & = & [\be_{\mathcal{D}^{\otimes k}} \otimes (\be_{\mathcal{D}} - d) + \be_{\mathcal{D}^{\otimes k}} \otimes d - d^{\otimes k} \otimes d] \\
& \le & [\be_{\mathcal{D}^{\otimes k}} \otimes (\be_{\mathcal{D}} - d)] + [(\be_{\mathcal{D}^{\otimes k}}  - d^{\otimes k}) \otimes \be_{\mathcal{D}}] \\
& \le & [(\be_{\mathcal{D}} - d) \otimes \be_{\mathcal{D}^{\otimes k}}] + k \cdot [(\be_{\mathcal{D}} - d) \otimes \be_{\mathcal{D}^{\otimes (k-1)}} \otimes \be_{\mathcal{D}}] \\
& = & (k+1) \cdot [(\be_{\mathcal{D}} - d) \otimes \be_{\mathcal{D}^{\otimes k}}]
\end{eqnarray*}
(using that $\mathcal{D}$ is strongly self-absorbing as well as  our induction hypothesis for the second inequality), so the assertion also  holds for $k+1$.
\end{nproof}
\en

\bn
\label{B}
The following is only a mild generalization of \cite[Lemma~1.3]{DadRor:ssa-projection}.

\begin{lms}
Let $\mathcal{D}$ be strongly self-absorbing and let $0 \le f \le g \le \be_{\mathcal{D}}$ be positive elements of $\mathcal{D}$ satisfying $\be_{\mathcal{D}} - g \neq 0$ and $fg = f$.

Then, there is $0\neq n \in \mathbb{N}$ such that 
\[
[f^{\otimes n}] \le [\be_{\mathcal{D}^{\otimes n}} - g^{\otimes n}] \mbox{ in } W(\mathcal{D}^{\otimes n}).
\]
\end{lms}

\begin{nproof}
Since $\mathcal{D}$ is simple and $\be_{\mathcal{D}} - g \neq 0$, there is $n \in \mathbb{N}$ such that 
\[
[f] \le n \cdot [\be_{\mathcal{D}} - g].
\]
Then, 
\begin{eqnarray*}
[f^{\otimes n}] & \le & n \cdot [(\be_{\mathcal{D}} - g) \otimes f \otimes \ldots \otimes f] \\
& = & [(\be_{\mathcal{D}}-g) \otimes f \otimes \ldots \otimes f] + \ldots + [f \otimes \ldots \otimes f \otimes (\be_{\mathcal{D}}-g)] \\
& = &  [(\be_{\mathcal{D}}-g) \otimes f \otimes \ldots \otimes f + \ldots + f \otimes \ldots \otimes f \otimes (\be_{\mathcal{D}}-g)] \\
& \le &  [(\be_{\mathcal{D}}-g) \otimes g \otimes \ldots \otimes g + \ldots + g \otimes \ldots \otimes g \otimes (\be_{\mathcal{D}}-g)] \\
& \le & [\be_{\mathcal{D}^{\otimes n}} - g^{\otimes n}],
\end{eqnarray*}
where for the first equality we have used that $\mathcal{D}$ is strongly self-absorbing, for the second equality we have used that the terms are pairwise orthogonal by our assumptions on $f$ and $g$, and the last inequality follows from Proposition~\ref{A}.
\end{nproof}
\en

\bn
\label{C}
The following is a version of \cite[Lemma~2.4]{DadRor:ssa-projection} for positive elements rather than projections.

\begin{lms}
Let $\mathcal{D}$ be strongly self-absorbing and let $0 \le f \le g \le \be_{\mathcal{D}}$ be positive elements satisfying $\be_{\mathcal{D}} - g \neq 0$ and $fg = f$; let $0 \neq d \in \mathcal{D}_{+}$.

Then, there is $0\neq m \in \mathbb{N}$ such that 
\[
[f^{\otimes m}] \le [d \otimes \be_{\mathcal{D}^{\otimes (m-1)}}] \mbox{ in } W(\mathcal{D}^{\otimes m}).
\]
\end{lms}

\begin{nproof}
By Lemma~\ref{B}, there is $0 \neq n \in \mathbb{N}$ such that
\[
[f^{\otimes n}] \le [\be_{\mathcal{D}^{\otimes n}} - g^{\otimes n}];
\]
since $f^{\otimes n} \perp \be_{\mathcal{D}^{\otimes n}} - g^{\otimes n}$, this implies 
\[
2 \cdot [f^{\otimes n}] \le [\be_{\mathcal{D}^{\otimes n}}].
\]
By an easy induction argument we then have
\[
2^{k} \cdot [f^{\otimes nk}] \le [\be_{\mathcal{D}^{\otimes nk}}]
\]
for any $k \in \mathbb{N}$. 

By simplicity of $\mathcal{D}$ and since $d$ is nonzero, there is $\bar{k} \in \mathbb{N}$ such that 
\[
[f] \le 2^{\bar{k}} \cdot [d].
\]
Set
\[
m:= n \bar{k} +1,
\]
then
\begin{eqnarray*}
[f^{\otimes m}] & \le & 2^{\bar{k}} \cdot [d \otimes f^{\otimes (m-1)}] \\
& = & 2^{\bar{k}} \cdot [d \otimes f^{\otimes n \bar{k}}] \\
& \le & [d \otimes \be_{\mathcal{D}^{\otimes n \bar{k}}}] \\
& = & [d \otimes \be_{\mathcal{D}^{\otimes (m-1)}}].
\end{eqnarray*} 
\end{nproof}
\en

\section{Large order zero maps}

\noindent
Below we establish the existence of nontrivial order zero maps from matrix algebras into strongly self-absorbing $C^{*}$-algebras, and we show certain systems of such maps give rise to order zero maps with small complements. We refer to  \cite{Win:cpr2} and \cite{WinZac:order-zero} for an introduction to order zero maps.

\bn
\label{G}
\begin{props}
Let $\mathcal{D}$ be strongly self-absorbing and $0 \neq d \in \mathcal{D}_{+}$.

Then, for any $0 \neq k \in \mathbb{N}$, there is a nonzero c.p.c.\ order zero map 
\[
\psi:M_{k} \to \overline{d\mathcal{D}d}.
\]
\end{props}

\begin{nproof}
Let us first prove the assertion in the case where $d = \be_{\mathcal{D}}$ and $k=2$. Since $\mathcal{D}$ is infinite dimensional, there are orthogonal positive normalized elements $e,f \in \mathcal{D}$. Since $\mathcal{D} \cong \mathcal{D} \otimes \mathcal{D}$ is strongly self-absorbing, there is a sequence of unitaries $(u_{n})_{n \in \mathbb{N}} \subset \mathcal{D} \otimes \mathcal{D}$ such that 
\[
u_{n} (e \otimes f) u_{n}^{*} \stackrel{n \to \infty}{\longrightarrow} f \otimes e;
\]
since $e \otimes f \perp f \otimes e$ this implies that there is a c.p.c.\ order zero map
\[
\textstyle
\bar{\sigma}: M_{2} \to \prod_{\mathbb{N}} \mathcal{D} \otimes \mathcal{D} / \bigoplus_{\mathbb{N}} \mathcal{D} \otimes \mathcal{D}
\]
given by 
\[
\bar{\sigma}(e_{11}) = e \otimes f,\; \bar{\sigma}(e_{22}) = f \otimes e, \; \bar{\sigma}(e_{21}) = \pi((u_{n}(e \otimes f))_{n \in \mathbb{N}})
\]
(cf.\ \cite{Win:cpr2}), where $\pi: \prod_{\mathbb{N}} \mathcal{D} \otimes \mathcal{D} \to \prod_{\mathbb{N}} \mathcal{D} \otimes \mathcal{D} / \bigoplus_{\mathbb{N}} \mathcal{D} \otimes \mathcal{D} $ denotes the quotient map. 

Since order zero maps with finite dimensional domains are semiprojective (cf.\ \cite{Win:cpr2}), $\bar{\sigma}$ has a c.p.c.\ order zero lift $M_{2} \to \prod_{\mathbb{N}}\mathcal{D} \otimes \mathcal{D}$, which in turn implies that there is a nonzero c.p.c.\ order zero map 
\[
\tilde{\sigma}:M_{2} \to \mathcal{D} \otimes \mathcal{D} \cong \mathcal{D}.
\]
Next, if $k= 2^{r}$ for some $r \in \mathbb{N}$, then 
\[
 M_{2^{r}} \cong (M_{2})^{\otimes r} \stackrel{\tilde{\sigma}^{\otimes r}}{\longrightarrow} \mathcal{D}^{\otimes r} \cong \mathcal{D}
\]
is a nonzero c.p.c.\ order zero map; for an arbitrary $k \in \mathbb{N}$, we may take $r$ large enough and restrict $\tilde{\sigma}^{\otimes r}$ to $M_{k} \subset M_{2^{r}}$ to obtain a nonzero c.p.c.\ order zero map
\[
\sigma:M_{k} \to \mathcal{D}.
\]
This settles the proposition for arbitrary $k$ and for $d = \be_{\mathcal{D}}$. Now if $d$ is an arbitrary nonzero positive element (which we may clearly assume to be normalized), we can define a c.p.c.\ map 
\[
\textstyle
\bar{\psi}:M_{k} \to \prod_{\mathbb{N}} \overline{d\mathcal{D} d} / \bigoplus_{\mathbb{N}} \overline{d\mathcal{D} d} \subset \prod_{\mathbb{N}}\mathcal{D} / \bigoplus_{\mathbb{N}} \mathcal{D}
\]
by setting
\[
\bar{\psi}(x):= \pi((d\sigma_{n}(x)d)_{n \in \mathbb{N}}) \mbox{ for } x \in M_{k},  
\]
where again $\pi: \prod_{\mathbb{N}} \overline{d\mathcal{D} d}\to \prod_{\mathbb{N}} \overline{d\mathcal{D} d} / \bigoplus_{\mathbb{N}} \overline{d\mathcal{D} d} $ denotes the quotient map and $\sigma_{n}:M_{k} \to \mathcal{D}$ is a sequence of c.p.c.\ maps lifting the c.p.c.\ order zero map 
\[
\textstyle
\mu \sigma: M_{k} \to (\prod_{\mathbb{N}} \mathcal{D} / \bigoplus_{\mathbb{N}} \mathcal{D}) \cap \mathcal{D}',
\]
with 
\[
\textstyle 
\mu : \mathcal{D} \to (\prod_{\mathbb{N}} \mathcal{D} / \bigoplus_{\mathbb{N}} \mathcal{D}) \cap \mathcal{D}'
\]
being a unital $*$-homomorphism as in \cite[Theorem~2.2]{TomsWin:ssa}. It is straightforward to check that $\bar{\psi}$ is nonzero and has order zero. Again by semiprojectivity of order zero maps, this implies the existence of a nonzero c.p.c.\ order zero map
\[
\psi:M_{k} \to \overline{d\mathcal{D}d}.
\]
\end{nproof}
\en

\bn
\label{H}
\begin{props}
Let $B$ be a unital $C^{*}$-algebra and $\varrho: M_{2} \to B$ a unital $*$-homomorphism. Define
\[
E:= \{f \in \mathcal{C}([0,1],B \otimes M_{2}) \mid f(0) \in B \otimes \be_{M_{2}}, \; f(1) \in \be_{B} \otimes M_{2}\}.
\]

Then, there is a unital $*$-homomorphism 
\[
\tilde{\varrho}:M_{2} \to E.
\]
\end{props}

\begin{nproof}
This follows from simply connecting the two embeddings $\varrho \otimes \be_{M_{2}}$ and $\be_{M_{2}} \otimes \id_{M_{2}}$ of $M_{2}$ into $\varrho(M_{2}) \otimes M_{2} \cong M_{2} \otimes M_{2}$ along the unit interval.
\end{nproof}
\en

\bn
\label{D}
\begin{lms}
Let $m \in \mathbb{N}$ and $A$ a unital $C^{*}$-algebra. Let
\[
\varphi_{1}, \ldots, \varphi_{m}:M_{2} \to A
\]
be c.p.c.\ order zero maps such that 
\[
\sum_{i=1}^{m} \varphi_{i}(\be_{M_{2}}) \le \be_{A} 
\]
and
\[
[\varphi_{i}(M_{2}), \varphi_{j}(M_{2})] = 0 \mbox{ if } i \neq j.
\]

Then, there is a c.p.c.\ order zero map 
\[
\bar{\varphi}:M_{2} \to C^{*}(\varphi_{i}(M_{2}) \mid i = 1, \ldots,m) \subset A
\]
such that
\[
\bar{\varphi}(\be_{M_{2}}) = \sum_{i=1}^{m} \varphi_{i}(\be_{M_{2}}).
\]
Moreover, if $d \in A_{+}$ satisfies $\varphi_{m}(e_{11})d=d$, we may assume that $\bar{\varphi}(e_{11})d = d$.
\end{lms}

\begin{nproof}
In the following, we write $C_{i}$, $i=1,\ldots,m$, for various copies of the $C^{*}$-algebra $\mathcal{C}_{0}((0,1],M_{2})$; these come equipped with c.p.c.\ order zero maps $\varrho_{i}:M_{2} \to C_{i}$ given by 
\[
\varrho_{i}(x)(t) = t \cdot x \mbox{ for } t \in (0,1] \mbox{ and } x \in M_{2}. 
\]
By \cite[Proposition~3.2(a)]{Win:cpr}, the c.p.c.\ order zero maps $\varphi_{i}:M_{2} \to A$ induce  unique $*$-homomorphisms $C_{i} \to A$ via $\varrho_{i}(x) \mapsto \varphi_{i}(x)$, for $x \in M_{2}$.

We now define a universal $C^{*}$-algebra
\[
\textstyle
B:= C^{*}(C_{i},\be \mid \sum_{l=1}^{m} \varrho_{l}(\be_{M_{2}}) \le \be, \; [C_{i},C_{j}] = 0 \mbox{ if } i \neq j \in \{1,\ldots,m\}).
\]
Then, $B$ is generated by the $\varrho_{i}(x)$, $i \in \{1, \ldots, m\}$ and $x \in M_{2}$; the assignment 
\[
\varrho_{i}(x) \mapsto \varphi_{i}(x) \mbox{ for } i \in \{1,\ldots,m\} \mbox{ and } x \in M_{2}
\] 
induces a unital $*$-homomorphism
\[
\pi:B \to C^{*}(\varphi_{i}(M_{2}),\be_{A} \mid i \in \{1,\ldots,m\}) \subset A
\]
satisfying 
\[
\sum_{l=1}^{m} \pi \varrho_{l}(\be_{M_{2}}) = \sum_{l=1}^{m} \varphi_{l}(\be_{M_{2}}). 
\]
Now if we find a c.p.c.\ order zero map 
\[
\bar{\varrho}:M_{2} \to B
\]
satisfying
\[
\bar{\varrho}(\be_{M_{2}}) = \sum_{l=1}^{m} \varrho_{l}(\be_{M_{2}}),
\]
then
\[
\bar{\varphi} := \pi \bar{\varrho}
\]
will have the desired properties, proving the first assertion of the lemma. We proceed to construct $\bar{\varrho}$.

For $k=1,\ldots,m$, let 
\[
\textstyle
J_{k}:= \mathcal{J}(\be - \sum_{l=k}^{m} \varrho_{l}(\be_{M_{2}})) \lhd B
\]
denote the ideal generated by $\be - \sum_{l=k}^{m} \varrho_{l}(\be_{M_{2}})$ in $B$; let 
\[
B_{k}:= B/J_{k}
\]
denote the quotient. We clearly have 
\[
J_{1} \subset J_{2} \subset \ldots \subset J_{m}
\]
and surjections
\[
B \stackrel{\pi_{1}}{\to} B_{1}   \stackrel{\pi_{2}}{\to} \ldots   \stackrel{\pi_{m}}{\to} B_{m}. 
\]
Observe that 
\[
\pi_{m} \circ \ldots \circ \pi_{1} \circ \varrho_{m} : M_{2} \to B_{m}
\]
is a unital surjective c.p.\ order zero map, hence a $*$-homomorphism by \cite[Proposition~3.2(b)]{Win:cpr}; therefore, 
\[
B_{m} \cong M_{2}.
\]
For $k=1,\ldots,m-1$, set
\begin{eqnarray}
\label{w1}
\lefteqn{E_{k}:= }\nonumber \\
&&  \{f \in \mathcal{C}([0,1 ], B_{k+1} \otimes M_{2}) \mid f(0) \in B_{k+1} \otimes \be_{M_{2}} , \; f(1) \in \be_{B_{k+1}} \otimes M_{2}  \} ;
\end{eqnarray} 
one easily checks that  the maps 
\begin{eqnarray*}
\sigma_{k} : B_{k} \to E_{k}  
\end{eqnarray*}
induced by 
\[
\pi_{k}  \ldots  \pi_{1}  \varrho_{i}(x) \mapsto
\left\{
\begin{array}{lcl}
(t & \mapsto & (1-t) \cdot \pi_{k+1}  \ldots  \pi_{1}  \varrho_{i}(x) \otimes \be_{M_{2}}) \\
&& \mbox{ for } i=k+1, \ldots, m \mbox{ and } x \in M_{2} \\
(t & \mapsto & t \cdot \be_{B_{k+1}} \otimes x) \\
&& \mbox{ for } i=k \mbox{ and } x \in M_{2}
\end{array}
\right.
\]
are well-defined  $*$-isomorphisms. Similarly, the map
\[
\sigma_{0}: B \to E_{0}:= \{f \in \mathcal{C}([0,1 ], B_{1}) \mid f(1) \in \mathbb{C} \cdot \be_{B_{1}} \}
\]
induced by
\begin{eqnarray*}
\varrho_{i}(x) & \mapsto & (t \mapsto (1-t) \cdot \pi_{1} \varrho_{i}(x)) \mbox{ for } i=1, \ldots,m \mbox{ and } x \in M_{2} , \\
\be_{B} & \mapsto & \be_{E_{0}}
\end{eqnarray*}
is a well-defined $*$-isomorphism; note that 
\[
\textstyle
\sigma_{0}(\sum_{l=1}^{m} \varrho_{l}(\be_{M_{2}}) ) = (t \mapsto (1-t) \cdot \be_{B_{1}}).
\]
By \eqref{w1} together with Proposition~\ref{H} and an easy induction argument, the unital $*$-homomorphism 
\[
\pi_{m} \ldots \pi_{1} \varrho_{m}:M_{2} \to B_{m}
\]
pulls back to a unital $*$-homomorphism 
\[
\tilde{\varrho}:M_{2} \to B_{1};
\]
This in turn  induces a c.p.c.\ order zero map 
\[
\check{\varrho}:M_{2} \to E_{0}
\] 
by
\[
\check{\varrho}(x) := (t \mapsto (1-t) \cdot \tilde{\varrho}(x));
\]
note that this map satisfies
\[
\check{\varrho}(\be_{M_{2}}) = (t \mapsto (1-t) \cdot \be_{B_{1}}).
\]
We now define a $*$-homomorphism
\[
\bar{\varrho}:= \sigma_{0}^{-1} \circ \check{\varrho}: M_{2} \to B;
\]
note that $\bar{\varrho}(\be_{M_{2}}) = \sum_{l=1}^{m} \varrho_{l}(\be_{M_{2}})$, whence $\bar{\varrho}$ is as desired.

For the second assertion of the lemma, note that $\bar{\varrho}$ and $\varrho_{m}$ agree modulo $J_{m}$. Therefore, $\bar{\varphi}= \pi \bar{\varrho}$ and $\varphi_{m} = \pi \varrho_{m}$ agree up to $\pi(J_{m})$. However, one checks that $\pi(J_{m}) \perp d$, whence $(\bar{\varphi}(x) - \varphi_{m}(x))d = 0$ for all $x \in M_{2}$. This implies $\bar{\varphi}(e_{11}) d = \varphi_{m}(e_{11})d = d$.
\end{nproof}
\en

\bn
\label{E}
\begin{props}
Let $\mathcal{D}$ be strongly self-absorbing, $0 \neq m \in \mathbb{N}$ and 
\[
\varphi_{0}:M_{2} \to \mathcal{D}
\]
a c.p.c.\ order zero map.

Then, there are c.p.c.\ order zero maps
\[
\varphi_{1}, \ldots, \varphi_{m}:M_{2} \to \mathcal{D}^{\otimes m}
\]
such that
\begin{enumerate}
\item[(i)] $\varphi_{1} = \varphi_{0} \otimes \be_{\mathcal{D}^{\otimes (m-1)}}$
\item[(ii)] $[\varphi_{i}(M_{2}),\varphi_{j}(M_{2})]= 0$ if $i \neq j$
\item[(iii)] $\be_{\mathcal{D}^{\otimes m}} - \sum_{i=1}^{m} \varphi_{i}(\be_{M_{2}}) = (\be_{\mathcal{D}} - \varphi_{0}(\be_{M_{2}}))^{\otimes m}$.
\end{enumerate}
\end{props}

\begin{nproof}
For $k \in \{1, \ldots,m\}$, define
\[
\varphi_{k}:= (\be_{\mathcal{D}}-\varphi_{0}(\be_{M_{2}}))^{\otimes (k-1)} \otimes \varphi_{0} \otimes \be_{\mathcal{D}^{\otimes (m-k)}},
\]
then the $\varphi_{k}$ obviously satisfy \ref{E}(i) and (ii).  

A simple induction argument shows that, for $k=1,\ldots,m$, 
\[
\be_{\mathcal{D}^{\otimes m}} - \sum_{i=1}^{k} \varphi_{i}(\be_{M_{2}}) = (\be_{\mathcal{D}}-\varphi_{0}(\be_{M_{2}}))^{\otimes k}  \otimes \be_{\mathcal{D}^{\otimes (m-k)}},
\]
which is \ref{E}(iii) when we take $k=m$.
\end{nproof}
\en

\section{$\mathcal{Z}$-stability}

\noindent
We now assemble the techniques of the preceding sections and a result from \cite{RorWin:Z-revisited} to prove our main result; we also derive some consequences.

\bn
\label{main}
\begin{thms}
Any strongly self-absorbing $C^{*}$-algebra  $\mathcal{D}$ absorbs the Jiang--Su algebra $\mathcal{Z}$ tensorially.
\end{thms}

\begin{nproof}
Let $k \in \mathbb{N}$. By Proposition~\ref{G}, there is a nonzero c.p.c.\ order zero map $\varphi:M_{2} \to \mathcal{D}$. Using functional calculus for order zero maps (cf.\  \cite{WinZac:order-zero}), we may assume that there is 
\[
0 \le d \le \varphi(e_{11})
\]
such that 
\[
d \neq 0 \mbox{ and }\varphi(e_{11})d = d.
\]
Note that 
\[
(\be_{\mathcal{D}}- \varphi(\be_{M_{2}})) (\be_{\mathcal{D}} - d) = \be_{\mathcal{D}} - \varphi(\be_{M_{2}}).
\]
By Proposition~\ref{G}, there is a nonzero c.p.c.\ order zero map 
\[
\psi:M_{k} \to \overline{d \mathcal{D}d}; 
\]
note that
\[
\varphi(e_{11}) \psi(x) = \psi(x) \mbox{ for } j=1, \ldots, k \mbox{ and } x \in M_{k}.
\]
Apply Lemma~\ref{C} (with $\mathcal{D}^{\otimes k}$, $\psi(e_{11})^{\otimes k}$, $(\be_{\mathcal{D}} - \varphi(\be_{M_{2}})) \otimes \be_{\mathcal{D}^{\otimes (k-1)}}$ and $(\be_{\mathcal{D}} - d) \otimes \be_{\mathcal{D}^{\otimes (k-1)}}$ in place of $\mathcal{D}$, $d$, $f$ and $g$, respectively) to obtain $0 \neq m \in \mathbb{N}$ such that 
\begin{equation}
\label{w2}
[((\be_{\mathcal{D}} - \varphi(\be_{M_{2}})) \otimes \be_{\mathcal{D}^{\otimes (k-1)}})^{\otimes m}] \le [\psi(e_{11})^{\otimes k} \otimes \be_{(\mathcal{D}^{\otimes k})^{\otimes  (m-1)}}]
\end{equation}
in $W((\mathcal{D}\otimes k)^{\otimes m})$. From Proposition~\ref{E} (with $\mathcal{D}^{\otimes k}$ in place of $\mathcal{D}$ and $\varphi_{0}:= \varphi \otimes \be_{\mathcal{D}^{\otimes (k-1)}}$) we obtain c.p.c.\ order zero maps 
\[
\varphi_{1}, \ldots, \varphi_{m}: M_{2} \to (\mathcal{D}^{\otimes k})^{\otimes m}
\]
satisfying \ref{E}(i), (ii) and (iii). By relabeling the $\varphi_{i}$ we may assume that actually $\varphi_{m} = \varphi_{0} \otimes \be_{(\mathcal{D}^{\otimes k})^{\otimes (m-1)}}$ in \ref{E}(i). 

From Lemma~\ref{D}, we obtain a c.p.c.\ order zero map 
\[
\bar{\varphi}: M_{2} \to C^{*}(\varphi_{i}(M_{2}) \mid i=1,\ldots,m) \subset (\mathcal{D}^{\otimes k})^{\otimes m}
\]
such that 
\[
\bar{\varphi}(\be_{M_{2}}) = \sum_{i=1}^{m} \varphi_{i}(\be_{M_{2}}).
\]
By the second assertion of Lemma~\ref{D} and since 
\begin{eqnarray*}
\varphi_{m}(e_{11}) (\psi(\be_{M_{k}}) \otimes \be_{\mathcal{D}^{\otimes (km-1)}})& = & (\varphi(e_{11})   \otimes \be_{\mathcal{D}^{\otimes (km-1)}})(\psi(\be_{M_{k}}) \otimes \be_{\mathcal{D}^{\otimes (km-1)}}) \\
& = & \psi(\be_{M_{k}}) \otimes \be_{\mathcal{D}^{\otimes (km-1)}},
\end{eqnarray*}
we may furthermore assume that 
\[
\bar{\varphi}(e_{11}) (\psi(\be_{M_{k}}) \otimes \be_{\mathcal{D}^{\otimes (km-1)}}) = \psi(\be_{M_{k}}) \otimes \be_{\mathcal{D}^{\otimes (km-1)}},
\]
which in turn yields 
\begin{equation}
\label{w4}
\psi(\be_{M_{k}}) \otimes \be_{\mathcal{D}^{\otimes (km-1)}} \le \bar{\varphi}(e_{11})
\end{equation}
since $\psi$ is contractive. Note that we have
\begin{eqnarray}
\label{w3}
[\be_{(\mathcal{D}^{\otimes k})^{\otimes m}} - \bar{\varphi}(\be_{M_{2}})] & \stackrel{\ref{E}\mathrm{(iii)}}{=} & [(\be_{\mathcal{D}^{\otimes k}} - \varphi_{0}(\be_{M_{2}}))^{\otimes m}]  \nonumber \\
& = &  [((\be_{\mathcal{D}} - \varphi(\be_{M_{2}})) \otimes \be_{\mathcal{D}^{\otimes (k-1)}})^{\otimes m}]  \nonumber \\
& \stackrel{\eqref{w2}}{\le} & [\psi(e_{11})^{\otimes k} \otimes \be_{(\mathcal{D}^{\otimes k})^{\otimes  (m-1)}}]
\end{eqnarray}
in $W((\mathcal{D}^{\otimes k})^{\otimes m})$. Define a c.p.c.\ order zero map
\[
\Phi: M_{2^{k}} \cong (M_{2})^{\otimes k} \to ((\mathcal{D}^{\otimes k})^{\otimes m})^{\otimes k}  \cong \mathcal{D}^{\otimes kmk}
\]
by
\[
\Phi:= \bar{\varphi}^{\otimes k}.
\]
We have
\begin{eqnarray*}
\lefteqn{ [\be_{((\mathcal{D}^{\otimes k})^{\otimes m})^{\otimes k}} - \Phi( \be_{(M_{2})^{\otimes k}})] } \\
& \stackrel{\ref{F}}{\le} & k \cdot [(\be_{(\mathcal{D}^{\otimes k})^{\otimes m}} - \bar{\varphi}(\be_{M_{2}})) \otimes \be_{((\mathcal{D}^{\otimes k})^{\otimes m})^{\otimes (k-1)}}]\\
& \stackrel{\eqref{w3}}{\le} & k \cdot [  \psi(e_{11})^{\otimes k} \otimes \be_{(\mathcal{D}^{\otimes k})^{\otimes  (m-1)}}    \otimes \be_{((\mathcal{D}^{\otimes k})^{\otimes m})^{\otimes (k-1)}}  ] \\
& \le & [ \psi(\be_{M_{k}})^{\otimes k} \otimes \be_{(\mathcal{D}^{\otimes (km - 1) })^{\otimes k}}] \\
& \stackrel{\eqref{w4}}{\le} & [\bar{\varphi}(e_{11})^{\otimes k}] \\
& = & [\Phi(e_{11})]
\end{eqnarray*}
in $W(((\mathcal{D}^{\otimes k})^{\otimes m})^{\otimes k})$. From \cite[Proposition~5.1]{RorWin:Z-revisited} we now see that there is a unital $*$-homomorphism
\[
\varrho: Z_{2^{k},2^{k}+1} \to \mathcal{D}^{\otimes kmk} \cong \mathcal{D}.
\]
Since $k$ was arbitrary, by \cite[Proposition~2.2]{TomsWin:ZASH} this implies that $\mathcal{D}$ is $\mathcal{Z}$-stable. 
\end{nproof}
\en

\bn
\label{initial-object-cor}
\begin{cors}
The Jiang--Su algebra is the uniquely determined (up to isomorphism) initial object in the category of strongly self-absorbing $C^{*}$-algebras (with unital $*$-homomorphisms).
\end{cors}

\begin{nproof}
By Theorem~\ref{main}, the Jiang--Su algebra does embed unitally into any strongly self-absorbing $C^{*}$-algebra, so it is an initial object. If $\mathcal{D}$ is another initial object, then $\mathcal{Z}$ and $\mathcal{D}$ embed unitally into one another, whence they are isomorphic by \cite[Proposition~5.12]{TomsWin:ssa}.  

Sometimes an object in a category is called initial only if there is a \emph{unique} morphism to any other object; this remains true in our setting if one takes approximate unitary equivalence classes of unital $*$-homomorphisms as morphisms.
\end{nproof}
\en

\bn
\label{K1-injective}
\begin{rems}
By \cite{Ror:Z-absorbing}, $\mathcal{Z}$-stable $C^{*}$-algebras are $K_{1}$-injective, whence $K_{1}$-injectivity is unnecessary in the hypotheses of the main results of \cite{TomsWin:ssa}, \cite{DadWinter:KK-D-topology}, \cite{DadWin:trivial-fields}, \cite{HirRorWin:D-stable}, \cite{HirWin:Rokhlin-ssa} and \cite{HirWin:permutations}. 
\end{rems}
\en


\begin{thebibliography}{10}

\bibitem{DadRor:ssa-projection}
Marius D\u{a}d\u{a}rlat and Mikael R{\o}rdam, \emph{Strongly self-absorbing
  {$C^*$}-algebras which contain a nontrivial projection}, arXiv preprint
  math.OA/0902.3886, 2009.

\bibitem{DadWinter:KK-D-topology}
Marius D\u{a}d\u{a}rlat and Wilhelm Winter, \emph{{On the $KK$-theory of
  strongly self-absorbing $C^{*}$-algebras}}, arXiv preprint math.OA/07040583,
  to appear in Math. Scand., 2007.

\bibitem{DadWin:trivial-fields}
Marius D\u{a}d\u{a}rlat and Wilhelm Winter, \emph{Trivialization of
  {$\mathcal{C} (X)$}-algebras with strongly self-absorbing fibres}, Bull. Soc.
  Math. France \textbf{136} (2008), no.~4, 575--606. \MR{MR2443037}

\bibitem{EllToms:regularity}
George~A. Elliott and Andrew~S. Toms, \emph{Regularity properties in the
  classification program for separable amenable {$C\sp *$}-algebras}, Bull.
  Amer. Math. Soc. (N.S.) \textbf{45} (2008), no.~2, 229--245. \MR{MR2383304}

\bibitem{HirRorWin:D-stable}
Ilan Hirshberg, Mikael R{\o}rdam, and Wilhelm Winter,
  \emph{{$\mathcal{C}_0(X)$}-algebras, stability and strongly self-absorbing
  {$C^*$}-algebras}, Math. Ann. \textbf{339} (2007), no.~3, 695--732.
  \MR{MR2336064 (2008j:46040)}

\bibitem{HirWin:Rokhlin-ssa}
Ilan Hirshberg and Wilhelm Winter, \emph{Rokhlin actions and self-absorbing
  {$C^*$}-algebras}, Pacific J. Math. \textbf{233} (2007), no.~1, 125--143.
  \MR{MR2366371}

\bibitem{HirWin:permutations}
\bysame, \emph{Permutations of strongly self-absorbing {$C^*$}-algebras},
  Internat. J. Math. \textbf{19} (2008), no.~9, 1137--1145. \MR{MR2458564}

\bibitem{JiaSu:Z}
Xinhui Jiang and Hongbing Su, \emph{On a simple unital projectionless
  {$C^*$}-algebra}, Amer. J. Math. \textbf{121} (1999), no.~2, 359--413.
  \MR{MR1680321 (2000a:46104)}

\bibitem{Ror:Z-absorbing}
Mikael R{\o}rdam, \emph{The stable and the real rank of
  {$\mathcal{Z}$}-absorbing {$C^*$}-algebras}, Internat. J. Math. \textbf{15}
  (2004), no.~10, 1065--1084. \MR{MR2106263 (2005k:46164)}

\bibitem{RorWin:Z-revisited}
Mikael R{\o}rdam and Wilhelm Winter, \emph{The {Jiang--Su} algebra revisited},
  arxiv preprint math.OA/0801.2259; to appear in J. Reine Angew. Math., 2008.

\bibitem{Toms:classproblem}
Andrew~S. Toms, \emph{On the classification problem for nuclear {$C\sp
  \ast$}-algebras}, Ann. of Math. (2) \textbf{167} (2008), no.~3, 1029--1044.
  \MR{MR2415391}

\bibitem{TomsWin:ssa}
Andrew~S. Toms and Wilhelm Winter, \emph{Strongly self-absorbing
  {$C^*$}-algebras}, Trans. Amer. Math. Soc. \textbf{359} (2007), no.~8,
  3999--4029 (electronic). \MR{MR2302521 (2008c:46086)}

\bibitem{TomsWin:ZASH}
\bysame, \emph{{$\mathcal{Z}$}-stable {ASH} algebras}, Canad. J. Math.
  \textbf{60} (2008), no.~3, 703--720. \MR{MR2414961}

\bibitem{Win:cpr}
Wilhelm Winter, \emph{Covering dimension for nuclear {$C^*$}-algebras}, J.
  Funct. Anal. \textbf{199} (2003), no.~2, 535--556. \MR{MR1971906
  (2004c:46134)}

\bibitem{Win:localizingEC}
\bysame, \emph{{Localizing the Elliott conjecture at strongly self-absorbing
  $C^{*}$-algebras}}, arXiv preprint math.OA/0708.0283v3, with an appendix by
  H.\ Lin, 2007.

\bibitem{Win:cpr2}
\bysame, \emph{{Covering dimension for nuclear $C^{*}$-algebras II}}, Trans.
  Amer. Math. Soc. \textbf{361} (2009), no.~8, 4143--4167.

\bibitem{WinZac:order-zero}
Wilhelm Winter and Joachim Zacharias, \emph{Completely positive maps of order
  zero}, arXiv preprint math.OA/0903.3290v1, 2009.

\end{thebibliography}

\providecommand{\bysame}{\leavevmode\hbox to3em{\hrulefill}\thinspace}
\providecommand{\MR}{\relax\ifhmode\unskip\space\fi MR }
\providecommand{\MRhref}[2]{%
  \href{http://www.ams.org/mathscinet-getitem?mr=#1}{#2}
}
\providecommand{\href}[2]{#2}

\end{document}